\documentclass[12pt]{article}

\usepackage{amsthm,amsmath,amssymb}

\usepackage{graphicx}
 \usepackage{epsfig}
  \usepackage{latexsym, epsfig, psfrag,eepic,colordvi,bm}
  \usepackage{graphics,color}
  \usepackage{graphics}
  \usepackage{graphicx}
  \usepackage{epsfig}
  \usepackage[all]{xy}

  \usepackage{amssymb}
  \usepackage{amsthm,amsmath}

  \usepackage{latexsym, epsfig, psfrag,eepic,colordvi,bm}
  \usepackage{graphics,color}
  \usepackage{amsmath,amsfonts,amssymb,amscd}
  \usepackage[all]{xy}
  \usepackage{epstopdf}

  \usepackage{amsthm,amsmath,amssymb}

  \usepackage{graphicx}
  \usepackage{amsthm,amsmath,amssymb}

  \usepackage{graphicx,epsfig}
\usepackage[colorlinks=true,citecolor=black,linkcolor=black,urlcolor=blue]{hyperref}

\usepackage[all]{xy}

\DeclareMathOperator{\Auto}{Aut}
\DeclareMathOperator{\Asyc}{Asyc}

\theoremstyle{plain}
\newtheorem{theorem}{Theorem}[section]

\newtheorem{lemma}{Lemma}[section]
\newtheorem{corollary}{Corollary}[section]
\newtheorem{proposition}{Proposition}[section]

\theoremstyle{definition}
\newtheorem{definition}[theorem]{Definition}

\theoremstyle{remark}

\date{}

\title{\bf The minimum number of chains in a noncrossing partition of a poset}

\author{Ricky X. F. Chen~\footnote{ORCID: 0000-0003-1061-3049}\\
	\small School of Mathematics, Hefei University of Technology\\[-0.8ex]
	\small Hefei, Anhui 230601, P.~R.~China\\[-0.8ex]
	\small\tt xiaofengchen@hfut.edu.cn
}

\begin{document}

\maketitle

\begin{abstract}
	The notion of noncrossing partitions of a partially ordered set (poset) is introduced here.
	When the poset in question is $[n]=\{1,2,\dots, n\}$ with the complete order of natural numbers,
	conventional noncrossing partitions arise. The minimum possible number of chains contained
	in a noncrossing partition of a poset clearly reflects the structural complexity of the poset. For the poset $[n]$, this number is just one. However,
	for a generic poset, it is a challenging task to determine the minimum number.
	Our main result in the paper is some characterization of this quantity.

  \bigskip\noindent \textbf{Keywords:}  partially ordered set; chain decomposition;
  noncrossing partition;
 linear extension; $132$-avoiding; descent

  \noindent\small Mathematics Subject Classifications: 06A07, 06A11, 05D05
\end{abstract}

\section{Introduction}

Partially ordered sets are well studied objects in discrete mathematics and
we will basically follow the notation in Stanley~\cite{stan1}.
A partially ordered set (poset) is a set $P$ with a binary relation `$\leq$' among the elements in $P$, where
the binary relation satisfies reflexivity, antisymmetry and transitivity. The poset will
be denoted by $(P,\leq)$ or $P$ for short.
For simplicity, all posets discussed in this paper are assumed to be finite.

If two elements $x$ and $y$ in $P$ satisfy $x\leq y$, we say $x$ and $y$ are comparable. We write $x<y$ if $x\leq y$ but $x\neq y$.
A chain of $P$ is a subset of elements such that any
two elements there are comparable, while an antichain
is a subset where any two elements are not comparable.
A chain decomposition of $P$ is a family of disjoint chains $\{C_1,C_2,\ldots, C_k\}$
such that $\bigcup_{i=1}^k C_i =P$. Let
$Min(P)$ denote the minimum number of chains that are contained in a chain
decomposition of $P$, and
let $Anti(P)$ denote the maximum number of elements that can be contained in an antichain of $P$.
These quantities reflect the structural complexity of the posets in question.
For instance, if there is a complete order in $P$, then $Min(P)=1$,
and if there is no order at all, $Min(P)=|P|$ (i.e., the number of elements in $P$).
The celebrated Dilworth's theorem~\cite{dil} states that $Min(P)=Anti(P)$ for all finite $P$.

Chain decompositions with various constraints
 have been studied, e.g., symmetric chain decomposition~\cite{greklei}, canonical
 symmetric chain decomposition~\cite{tim}, etc.,
which reflect the structural properties
and complexity of posets from different angles.
Here we introduce a new chain decomposition which can be viewed as
a generalization of the ubiquitous object in combinatorics, i.e., noncrossing partitions (e.g., see Armstrong~\cite{arm} and
Simion~\cite{simion}).
As such, we call the new decompositions noncrossing partitions of posets.
Specifically, a noncrossing partition of the set $[n]=\{1,2,\ldots, n\}$
is merely a noncrossing partition of the poset $[n]$ with the natural order.
Note that $[n]$ itself is a noncrossing partition of $[n]$.
That is, the minimum number of chains contained in a noncrossing partition is simplely one in this case.
However, determining the minimum number of chains in a noncrossing partition for a general poset is a challenging
task.

Our main goal of this note is to provide some characterization of the minimum possible number of chains
contained in a noncrossing partition of a generic poset.

\section{Noncrossing partitions of posets}

Recall a noncrossing partition (see~\cite{arm,simion}) of the set $[n]$ is a partition
of $[n]$ into $k\geq 1$ blocks $B_1, B_2, \ldots, B_k$
such that there do not exist elements $a, b\in B_i$ and $c, d \in B_j$ ($i\neq j$)
such that $a<c<b<d$.
{For example, for $n=5$, $B_1=\{1,5\}$ and $B_2=\{2,3,4\}$
give a noncrossing partition of $[5]$ into two blocks, while $B_1=\{1,3,5\}$ and
$B_2=\{2,4\}$ do not give a noncrossing partition.}
Evidently, the definition depends on the natural order on $[n]$.
Regarding $[n]$ as a poset, $B_i$ is just a chain and
a partition is just a chain decomposition.
What if we replace $[n]$ with an arbitrary poset?

{Before we proceed, it will be convenient to represent a poset by introducing its Hasse diagram.
For two elements $x$ and $y$ in a poset $P$, if $x< y$ and there does not exist $z$ such that
$x< z <y$, then we say $y$ covers $x$. The Hasse diagram of $P$ is the graph with the elements of $P$
as the vertices, and with the covering relation giving the edges, and if $y$ covers $x$ then $y$ is drawn ``above" $x$ (with an edge
between $x$ and $y$).
Note the whole partial relation can be derived by applying the transitivity based on the Hasse diagram.}

\begin{definition}
	A chain decomposition of a poset $P$, $\{C_1, C_2, \ldots, C_k\}$, is called a noncrossing partition of $P$, if there do not exist elements $a, b\in C_i$ and $c, d \in C_j$ ($i\neq j$)
	such that $a<c<b<d$ in $P$.
\end{definition}

\begin{figure}[!htb]
	\centering
	\includegraphics[scale=.66]{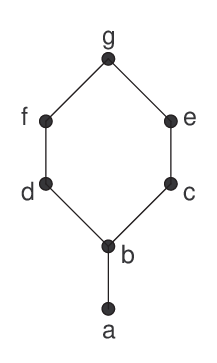}
	\caption{A poset of $7$ elements represented by its Hasse diagram.}\label{fig:hass-1}
\end{figure}

{
	For example, for the poset in Figure~\ref{fig:hass-1}, $\{\{a,b,c,e\}, \{d,f,g\}\}$ is a noncrossing partition,
	while $\{\{a,c,e\}, \{b,d,f,g\}\}$ is not since $a<b<c<g$.
}
We denote by $Min_{nc}(P)$ the minimum number of chains contained in a noncrossing
partition of $P$. {\color{blue}For $P$ in Figure~\ref{fig:hass-1}, $Min_{nc}(P)=2$.}
Clearly, $Min_{nc}(P)=1$ if and only if $P\sim [n]$.
However,
it is not easy to exactly determine this number for a generic poset.
Nevertheless, by relating noncrossing partitions to other notion,
we are able to prove some bounds.

\begin{definition}
	Let $(P,\leq)$ be a poset. A homogeneous chain decomposition (HCD) $\mathcal{C}$ of $P$ is a collection of mutually disjoint
	chains $C_1,C_2,\ldots, C_n$ such that $\bigcup_i C_i=P$, and if $s_i\in C_i$ and $s_j\in C_j$ are
	comparable, then all elements in $C_i$ and $C_j$ are pairwise comparable.
\end{definition}
{In the example of Figure~\ref{fig:hass-1}, $\{\{a,b,g\}, \{c, e\}, \{d,f\}\}$
gives an HCD.}
When all elements in two chains $C_i$ and $C_j$ are pairwise comparable, i.e.,
$C_i \bigcup C_j$ is a chain, we say $C_i$ and
$C_j$ are comparable for short. We also write
$C=(\xi_1<\xi_2<\dots<\xi_s)$ as a shorthand of that $C$ is the chain $\{\xi_1,\xi_2,\ldots, \xi_s\}$ and
$\xi_1<\xi_2<\dots<\xi_s$.

Denote by $|\mathcal{C}|$ the number of chains contained in $\mathcal{C}$.
Let $Min_h(P)=\min_{\mathcal{C}} |\mathcal{C}|$, where the minimization is over all HCDs of $P$.
An HCD of $P$ containing exactly $Min_h(P)$ chains is called a minimal homogeneous chain decomposition (MHCD) of $P$.
It turns out there is only one such a decomposition.
{For instance, for $P$ in Figure~\ref{fig:hass-1}, $Min_h(P)=3$ and $\{\{a,b,g\}, \{c, e\}, \{d,f\}\}$
	is actually the only MHCD.}

\begin{proposition}\label{2thm1}
	For any poset $P$, there exists a unique MHCD of $P$.
\end{proposition}
\proof Let $\mathcal{C}=\{C_1,C_2,\ldots, C_m\}$ be a MHCD of $P$.
If $|C_i|=1$ for all $1\leq i \leq m$, there is nothing to prove.
Thus, we assume that there exists at least one $i$ such that $|C_i|>1$.
Let $\mathcal{C'}=\{C'_1,C'_2,\ldots, C'_m\}$ be a different MHCD of $P$.
First, there exists $k$ and $s_{k1},\: s_{k2}\in C_k$ such that $s_{k1}\in C'_{j1}$, $s_{k2}\in C'_{j2}$
and $C'_{j1}\neq C'_{j2}$.
Otherwise, it is not hard to argue $\mathcal{C}=\mathcal{C'}$.
Next, since $s_{k1}$ and $s_{k2}$ are comparable, $C'_{j1}$ and $C'_{j2}$ are comparable.
Thus, $C^*=C'_{j1} \bigcup C'_{j2}$ is a chain of $P$.

We claim $(\mathcal{C'}\setminus \{C'_{j1},C'_{j2}\}) \bigcup \{C^*\}$ is an HCD of $P$.
For any $j\notin \{j1,j2\}$, $C'_j$ is either comparable to $C'_{j1}$ or not comparable to $C'_{j1}$.
For the former case, there exists $s_j\in C'_j$ comparable to $s_{k1}$. Since $s_{k1}$ and $s_{k2}$ come from
the same chain $C_k$, regardless of whether $s_j\in C_k$, $s_{k2}$ must be comparable to $s_j$ as well.
Hence, $C'_j$ is also comparable to $C'_{j2}$ so that $C'_j$ is comparable to $C^*$.
For the latter case, we can analogously show $C'_j$ is not comparable to $C^*$.
Thus, the claim holds.
However, this contradicts the assumption that $\mathcal{C'}$ is minimum.
Hence, $\mathcal{C}$ is the unique MHCD of $P$.\qed

HCDs
were first introduced in Chen and Reidys~\cite{chr-1} in the context of studying
the interaction between incidence algebras of posets and linear sequential dynamical systems,
where in particular, it was shown that the M\"{o}bius function of
any poset can be efficiently computed via a sequential dynamical system and a cut theorem concerning
HCDs of posets holds.

Another notion that we need is a generalization of $132$-avoiding permutations, another
ubiquitous object in combinatorics and computer science.

\begin{definition}
A permutation $\pi=\pi_1 \pi_2 \cdots \pi_n$ of the elements of $P$
is called $132$-avoiding if no three-element subsequence $\pi_{i_1} \pi_{i_2} \pi_{i_3}$ in $\pi$
satisfies $i_1<i_2<i_3$ while $\pi_{i_1}< \pi_{i_3} < \pi_{i_2}$ in $P$.
\end{definition}

In the case of $P=[n]$,
conventional $132$-avoiding permutations arise.
{\color{blue}For example, when $P=[5]$, $53241$ is a $132$-avoiding permutation, while $2 1 4 5 3$ is not.
Because in the latter, we realized that the subsequences $243$, $253$, $143$, $153$ all violate the definition.}
A linear extension of $P$ is a permutation $e=e_1e_2\cdots e_n$ of $P$-elements
such that $e_i<e_j$ implies $i<j$.
{\color{blue} For example, for the poset $P$ in Figure~\ref{fig:hass-1}, $abcdefg$ and $abcedfg$ are linear extensions.}
There are more than one linear extension unless $P\sim [n]$.

\begin{definition}
	Let $e=e_1e_2\cdots e_n$ be a linear extension of $P$.
A permutation $\pi=\pi_1 \pi_2 \cdots \pi_n$ of the elements of $P$
is called $132^e$-avoiding (i.e., $132$-avoiding with respect to $e$) if there does not exist a subsequence $\pi_{i_1} \pi_{i_2} \pi_{i_3}= e_{j_1} e_{j_2} e_{j_3}$
such that $i_1<i_2<i_3$ and $j_1< j_3 <j_2$.
\end{definition}

{\color{blue} For example, $gedfbac$ is $132$-avoiding w.r.t.~the linear extension $abcdefg$ of $P$ in Figure~\ref{fig:hass-1},
while $abcedfg$ is not due to the appearance of the subsequence $ced$.}
It is easily seen that a $132^e$-avoiding permutation is a $132$-avoiding permutation of $P$.
Given a permutation $\pi=\pi_1 \pi_2 \cdots \pi_n$ of $P$, $i$ is called a p-descent of $\pi$
if $\pi_i>\pi_{i+1}$ or $\pi_i$ is not comparable with $\pi_{i+1}$ in $P$ or $i=n$.
The number of p-descents in $\pi$ is denoted by $d_P(\pi)$.
{\color{blue} For $P$ in Figure~\ref{fig:hass-1} and $\pi=gedfbac$, it can be checked that $d_P(\pi)=5$, i.e., $i=1,2,4,5,7$.}
Let 
\begin{align*}
	Min_d(P) &=\min\{d_P(\pi): \mbox{$\pi$ is a $132$-avoiding permutation of $P$}\},\\
	Min_d^e(P) &=\min\{d_P(\pi): \mbox{$\pi$ is a $132^e$-avoiding permutation of $P$}\}.
\end{align*}
Now we are in a position to present our main result.
\begin{theorem}[Main theorem]\label{thm:nc-better}
	For any poset $P$, there exists a linear extension $e$ of $P$ such that
	\begin{align}
	Min_{nc}(P) \leq Min_d(P) \leq Min_d^e(P) \leq Min_h(P).
	\end{align}
Moreover, all inequalities are sharp.
\end{theorem}

\begin{figure}[!htb]
	\centering
	\includegraphics[scale=.66]{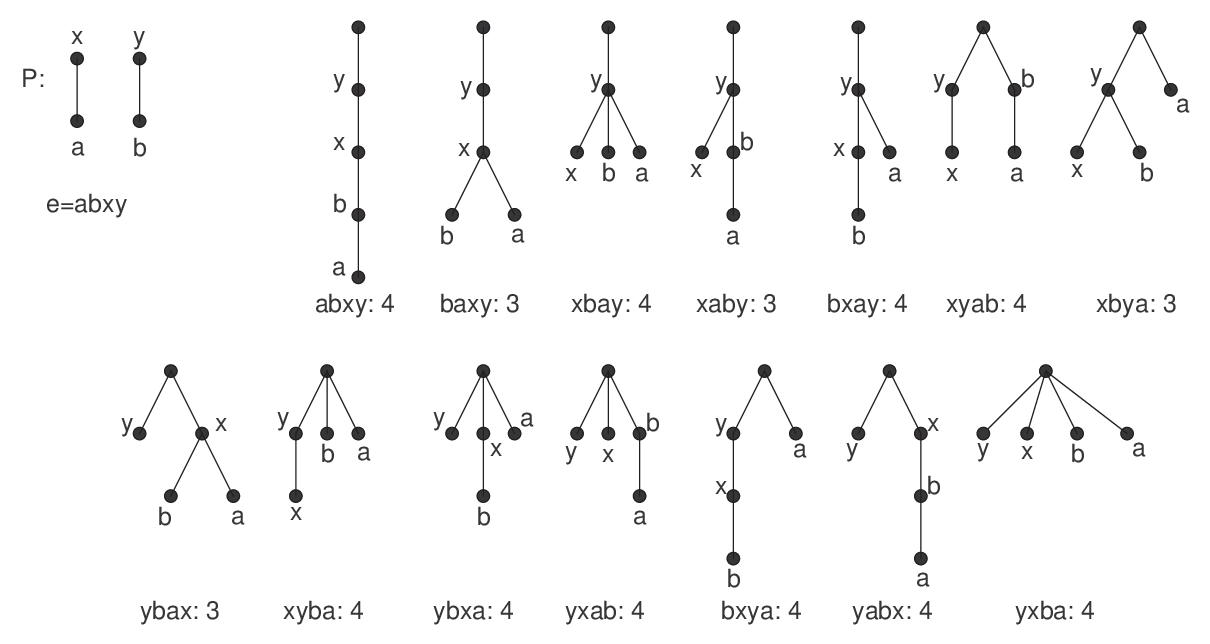}
	\caption{A poset $P$ with a linear extension $e$ such that $Min_d^e(P)> Min_h(P)$.}\label{fig:counterex}
\end{figure}

We remark that the rightmost inequality is not necessarily true for an arbitrary linear extension.
{For example, for the poset $P$ in Figure~\ref{fig:counterex}, it is easy to see $Min_h(P)=2$.
However, for its linear extension $e=abxy$, there are $14$ $132^e$-avoiding permutations all of which
have at least three p-descents. In Figure~\ref{fig:counterex}, the number of descents of a $132^e$-avoiding permutation
is written right after the $132^e$-avoiding permutation. For instance, ``$baxy: 3$" means that $baxy$ has three p-descents. }
Moreover, the reason that we are interested in $Min_d^e(P)$ is as follows:
while it may be hard to generate all $132$-avoiding permutations of $P$ to compute $Min_d(P)$,
it is easy to generate all $132^e$-avoiding permutations for any linear extension $e$ of $P$
as we shall see it is essentially generating all plane trees.
A proof of the above theorem follows from a series of properties that we are about to present.

Assume $\pi=\pi_1\pi_2 \cdots \pi_n$ is a permutation of a poset $P$.
Read $\pi$ from left to right and collect these elements between two consecutive p-descents,
excluding the first one and including the second.
By definition of p-descents, these elements comprise a chain.
In this way, all p-descents of $\pi$ induce a chain decomposition of $P$.

\begin{proposition}\label{prop:p-descent}
	Let $\pi$ be a $132$-avoiding permutation of a poset $P$.
	Then, the induced chain decomposition by the p-descents of $\pi$ is a noncrossing partition of $P$.
\end{proposition}
\proof 
If not, without loss of generality, suppose $\pi_1, \pi_2$ from the first induced chain and $\pi_3, \pi_4$
from the second induced chain cross, i.e.,
$\pi_1 < \pi_3 < \pi_2 <\pi_4$ or $\pi_3 < \pi_1 < \pi_4 < \pi_2$.
Obviously, either case implies a $132$ pattern in $\pi$, a contradiction whence the proposition. \qed

As a result, we immediately have $Min_{nc}(P)\leq Min_d(P) \leq Min_d^e(P)$ for any linear extension $e$ of $P$.
If otherwise explicitly stated, we assume the following notation in the rest of the section.
Let $\mathcal{C}=\{C_1,C_2,\ldots, C_k\}$ be the MHCD of $P$, where 
$$
C_i=(s_{i1}<s_{i2}<\dots<s_{im_i}), \quad
\sum_{i=1}^k m_i =n.
$$
 
\begin{lemma}\label{keylem}
	If $C_i$ and $C_j$ are comparable, then there exists $0\leq l \leq m_i$
	such that
	$$
	s_{i1}<s_{i2}<\cdots <s_{il}<s_{j1}<s_{j2}<\cdots <s_{jm_j}<s_{i(l+1)}<s_{i(l+2)}<\cdots < s_{im_i}.
	$$
\end{lemma}
\proof In order to prove the lemma, it suffices to show that there does not exist $0<l_1<m_i$ and $0<l_2<m_j$ such
that
$$
s_{il_2}<s_{j{l_1}}<s_{i(l_2+1)}< s_{j(l_1+1)} .
$$
Assume by contradiction that such $l_1$ and $l_2$ exist.
For any other chain $C_k$, if $C_k$ is comparable to $C_i$ and $x\in C_k$,
then either $x<s_{i(l_2+1)}$ or $x>s_{i(l_2+1)}$.
In any case, we conclude that an element in $C_j$ is comparable to $x$
whence $C_j$ and $C_k$ are comparable.
By similar analysis, we can conclude that if $C_k$ is not comparable to $C_i$,
then $C_k$ is not comparable to $C_j$ either.
Therefore, $(\mathcal{C}\setminus \{C_i,C_j\})\bigcup \{C_i\cup C_j\}$ is
an HCD of $P$.
This contradicts the assumption that $\mathcal{C}$ is the minimum and the lemma follows. \qed

Consider the relation $\leq_b$ on $\mathcal{C}$ that
$C_i \leq_b C_j$ if there exist elements $x,z\in C_j$ and $y\in C_i$
such that $x<y<z$ or $\min(C_i)>\max(C_j)$.
As for the first case, we say $C_j$ wrap around $C_i$ or $C_i$
can be wrapped around by $C_j$.
In view of Lemma~\ref{keylem}, we leave it to the reader to verify that
$(\mathcal{C},\leq_b)$ is a well-defined poset.

\begin{proposition}\label{prop:132}
	Suppose $C_1C_2\cdots C_k$ is a linear extension of $(\mathcal{C}, \leq_b)$.
	Then the following permutation $\pi$ is $132$-avoiding and has $k$ p-descents:
	$$
	\pi=	s_{11}s_{12}\cdots s_{1m_1} s_{21}\cdots s_{2m_2}\cdots s_{k1}\cdots s_{km_k}.
	$$
\end{proposition}
\proof 
By definition, it is easy to see there are exactly $k$ p-descents in $\pi$.
We prove the rest by contradiction.
Suppose $ \pi_{l_1} \pi_{l_2} \pi_{l_3}$ is a $132$ pattern in $\pi$. Since each $C_i$ appears as an increasing chain in $\pi$,
we have only two possible cases:
\begin{itemize}
\item $\pi_{l_1}, \pi_{l_2} \in C_i$, $\pi_{l_3} \in C_j$, and $i<j$;
\item $\pi_{l_1}\in C_i, \pi_{l_2} \in C_j$, $\pi_{l_3} \in C_k$, and $i<j<k$.
\end{itemize}
The first case cannot happen because the condition implies that $C_j <_b C_i$ in the light of Lemma~\ref{keylem},
contradicting the assumption of the proposition.
Next suppose the second case occurs.
First, $\pi_{l_1}< \pi_{l_2}$ and $C_i <_b C_j$ imply that
$\pi_{l_1} > x_j$ for some $x_j \in C_j$,
i.e., $C_j$ wrap around $C_i$.
Analogously, $C_k$ wrap around $C_i$.
Secondly, $\pi_{l_2} > \pi_{l_3}$ and $C_j <_b C_k$ imply that either
$\min(C_j)>\max(C_k)$ or $C_k$ wrap around $C_j$.
Since both $C_j$ and $C_k$ can wrap around $C_i$,
the former is absurd.
On the other hand, that $C_k$ wrap around $C_j$ while $C_j$ wrap around $C_i$
makes it impossible to have a $132$ pattern $ \pi_{l_1} \pi_{l_2} \pi_{l_3}$ such
that $\pi_{l_1}\in C_i, \pi_{l_2} \in C_j$, $\pi_{l_3} \in C_k$.
Hence, no $132$ patterns exist in $\pi$, completing the proof. \qed

From Lemma~\ref{keylem} and Proposition~\ref{prop:132}, we conclude 
$$
Min_{nc}(P) \leq Min_d(P) \leq Min_h(P).
$$
But we cannot conclude $Min_d^e(P)\leq Min_h(P)$ for an arbitrary linear extension $e$. 

\begin{figure}[!htb]
	\centering
	\includegraphics[scale=.66]{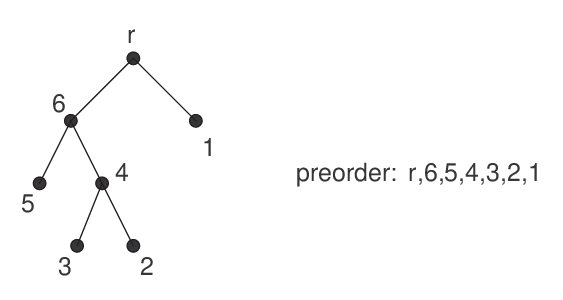}
	\caption{A labelled plane tree and the preorder of its vertices.}\label{3fig4}
\end{figure}

We proceed with further analysis below,
where on the way we need to make use of plane trees.
A {plane tree} $T$ can be recursively defined as an unlabeled tree with one distinguished vertex called the {root} of $T$, where the unlabeled trees obtained by deleting the root as well as its adjacent edges from $T$ are linearly ordered, and they are plane trees with the vertices adjacent to the root of $T$ in $T$ as their respective roots.
These subtrees are pictured as locating below the root and appearing left to right.
 A non-root vertex without any child is called a {leaf}, and an {internal vertex} otherwise.
 The root is always treated as internal.
A labelled plane tree is a plane tree where vertices carry mutually distinct labels from a certain set of labels.
The preorder of the vertices in a labelled plane tree $T$ is the sequence obtained by
travelling $T$ in a left-to-right depth-first manner and recording the label of a vertex when it is first visited.
See Figure~\ref{3fig4} for an example.

There is a bijection between plane trees and conventional $132$-avoiding permutations given by Jani and Rieper~\cite{jani-rieper}.
The following is how it works.
Let $T$ be a plane tree of $n$ edges. We use a preorder traversal of $T$ to label the non-root vertices with the integers $n, n-1,\ldots, 1$. As such, the first vertex visited gets the label $n$ and the last receives 1. A permutation written as a word is next obtained by reading the labelled tree in postorder, that is, traverse the tree from left to right and record the label of a vertex when it is last visited.
It was shown~\cite{jani-rieper} that the obtained permutation is $132$-avoiding on $[n]$.
As an example, for the tree in Figure~\ref{3fig4}, the obtained $132$-avoiding permutation is $532461$.

The reverse from a $132$-avoiding permutation to a plane tree was not explicitly presented in Jani and Rieper~\cite{jani-rieper}.
A reverse procedure was proposed in~\cite{chen2} and is presented here
for later use.
Let $\pi$ be a $132$-avoiding permutation on $[n]$.
Suppose the (increasing) chains induced by the p-descents of $\pi$ from left to right are $\tau_1, \tau_2, \ldots, \tau_k$.
Start with $\tau_k$ and make it into a path with the maximum (i.e., rightmost) element in $\tau_k$ attaching to the root of the expected plane tree $T$.
For example, suppose $\pi=532461$. Then, $k=4$ and $\tau_k=1$, and the path will be the path from vertex $1$ to the root of the tree in
Figure~\ref{3fig4}.
After $\tau_i$ is ``integrated" into the (partial) tree, we find the minimum element $u$ in the
path from the leftmost leaf to the root in the current partial tree that is larger than the maximum element $x$ in $\tau_{i-1}$,
and attach the path induced by $\tau_{i-1}$ to the tree such that $u$ and $x$ are adjacent;
if no such a $u$ exists, we attach the path induced by $\tau_{i-1}$ to the root of the current tree.
Iteration of the procedure eventually yields a labelled plane tree.
(The vertex labels are uniquely determined by the underlying plane tree.)

In the forthcoming result, a straightforward generalization of $132^e$-avoiding permutations of $P$ from
a linear extention $e$ to an arbitrary permutation $e$ of $P$ will be used.

\begin{proposition}\label{prop:e}
	Suppose $C_1C_2\cdots C_k$ is
	a linear extension of the poset $(\mathcal{C},\leq_b)$.
	Then, there exists a labelled plane tree $T$ with non-root vertex labels from $P$ such that
	$\pi$ in Proposition~\ref{prop:132} is $132^e$-avoiding, where $e$ is the reverse
	of the preorder of the vertices other than the root of $T$.
\end{proposition}
\proof First, we use the chains $C_i$ to build a labelled plane tree following
the same procedure from the induced chains of $132$-avoiding permutations to plane trees described above.
We then claim the obtained tree is the desired $T$.
To see this, one has to realize that the Jani-Rieper bijection essentially
encodes the relation among the non-root vertex sequences from the preorder, postorder and the reverse of the preorder.
In a word, the postorder is $132$-avoiding with respect to the reverse of the preorder.
{\color{blue} Actually, this is how we constructed all $132^e$-avoiding permutations in Figure~\ref{fig:counterex}.}
In particular, when the preorder is $n, n-1, \ldots, 1$, the postorder gives a conventional
$132$-avoiding permutation on $[n]$.
 The rest is clear, completing the proof. \qed

It remains to prove that there exists a linear extension of $(\mathcal{C},\leq_b)$ of which the corresponding $e$
is in fact a linear extension of $P$.
We need one more important lemma to that end.

\begin{lemma}\label{lem:maximal}
	Suppose $\{C_{i_1}, C_{i_2}, \ldots, C_{i_{k'}}\}$ is a subposet of $(\mathcal{C},\leq_b)$,
	and $C_{i_1},C_{i_2},\ldots, C_{i_m}$ are the maximal elements of the subposet.
	Then, any $C_{i_j}$ for $m+1\leq j \leq k'$ satisfies either one of the cases:
	\begin{itemize}
		\item[$(1)$] for at least one $t$ ($1\leq t \leq m$), $\min(C_{i_j})>\max(C_{i_t})$;
\item[$(2)$] for a unique $t$ ($1\leq t \leq m$), $C_{i_t}$ wrap around $C_{i_j}$.
	\end{itemize}
	In addition, two case $(2)$ elements wrapped around by distinct maximal elements are not comparable,
	while a case $(1)$ element is smaller than a case $(2)$ element if comparable and
	the minimal ($P$-element) of the former is greater than the maximal of the latter.
\end{lemma}
\proof For any $m+1\leq j \leq k'$, by assumption, $C_{i_j}$ is smaller than at least one maximal element.
We first show that $C_{i_j}$ cannot satisfy both cases.
Suppose $\min(C_{i_j})>\max(C_{i_t})$ for some $1\leq t \leq m$.
If $C_{i_j}$ can be wrapped around by another maximal element $C_{i_{t'}}$,
then it is easy to see that $C_{i_t}$ and $C_{i_{t'}}$ are comparable, a contradiction.
Analogously, an element satisfying $(2)$ cannot satisfy $(1)$ at the same time.
Moreover, an element cannot be wrapped around by more than one maximal element.

If two case $(2)$ elements wrapped around by distinct maximal elements are comparable, either
the minimal $P$-element (i.e., element in $P$) of one is greater than the maximal $P$-element of the other or one wrap around the other.
Either case implies the two involved distinct maximal elements are comparable, contradicting
the maximality.
The remaining statement can be similarly verified, and the proof follows.
\qed

\begin{proposition}
	There exists a linear extension of $(\mathcal{C},\leq_b)$, still denoted by $C_1C_2\cdots C_k$, of which the corresponding $e$ referred to in Proposition~\ref{prop:e}
	is a linear extension of $P$.
\end{proposition}
\proof Our strategy here is to construct a linear extension of $(\mathcal{C},\leq_b)$
first and then argue the corresponding $e$ is a linear extension of $P$.

\begin{figure}[!htb]
	\centering
	\includegraphics[scale=.66]{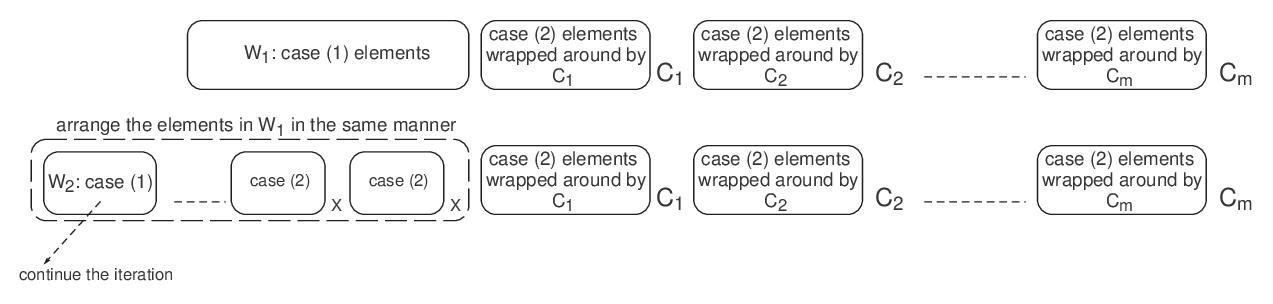}
	\caption{Construct a linear extension of $(\mathcal{C},\leq_b)$.}\label{fig:e-constr}
\end{figure}

{\bf Construct a linear extension of $(\mathcal{C},\leq_b)$.}
First, we apply the following procedure. 
\begin{itemize}
	\item[(i)] Initialize $j=0$, $W_0=(\mathcal{C},\leq_b)$ and $F=W_0$;
\item[(ii)] Arrange the maximal elements of $F$ on a line in an arbitrary way;
\item[(iii)] Put each of those case $(2)$ elements w.r.t.~$F$ right before the maximal element in $F$ that wrap around it (and
after the preceeding maximal element) and order those right before the same maximal element later.
Next, denote the set of case $(1)$ elements w.r.t.~$F$ by $W_{j+1}$ and put them in front of the current ``partially linearized" sequence
in an arbitrary way. Update $F$ to $W_{j+1}$ and $j$ to $j+1$;
\item[(iv)] Iterate (ii) and (iii) until $F$ is an empty set.
\end{itemize}
At this point, all elements of $(\mathcal{C}, \leq_b)$ are in a sense grouped into disjoint ordered groups.
The involved maximal $\mathcal{C}$-elements (w.r.t.~a certain iteration) serve as a kind of group markers.
(The group marker of a group is on the right-hand side.)
See Figure~\ref{fig:e-constr} for an illustration,
where $C_m$ and the case (2) elements wrapped around by $C_m$
give an example of a group.
The groups obtained so far will be referred to as type I groups.
By construction, the maximum $P$-element contained in a group marker is larger than (in terms of $(P,\leq)$) all other
$P$-elements contained in the chains in the same group.
Moreover, in view of Lemma~\ref{lem:maximal}, any $\mathcal{C}$-element in a left group is smaller than
any $\mathcal{C}$-element in a right group if comparable, not violating the current sequence to possibly become a linear extension of $(\mathcal{C},\leq_b)$.

Iteratively apply the above procedure to each type I group with the group marker excluded and each of
those newly generated groups (excluding the group markers) in the process until each non-empty group contains a single element.
It is a kind of successive ``linearization".
Eventually, we obtain a linear extension of $(\mathcal{C}, \leq_b)$.

Assume the resulting linear extension is $C_1 C_2\cdots C_k$,
and its corresponding tree from the reverse procedure of the Jani-Rieper bijection is $T$.
Note that in terms of Figure~\ref{fig:e-constr}, $C_k$ here is actually $C_m$, i.e., the
rightmost chain.
We next show that {\bf the reverse $e$ of the preorder of the vertices other than the root of $T$
is a linear extension of $P$,}
which is equivalent to showing that for any two entries
in the preorder, the left one is greater than the right one if comparable in $P$. 
To this end, for any vertex
$u$ in $T$, consider the subtree $T_u$ with $u$ as the root. It suffices to show: (i) $u$ is greater than any of its
descendants (in terms of the labels in $P$) where the root of $T$ is assumed to be
an artificial maximum element added into $P$; (ii) any vertex in a left
subtree of $T_u$ is greater than any vertex in a right subtree of $T_u$
if comparable.

Suppose $u$ is the root of $T$. Then, $T_u=T$.
According to the construction of the linear extension $C_1 C_2 \cdots C_k$
and the tree $T$, $P$-elements contained in chains belonging to distinct type I groups (including respective group markers)
are contained in distinct subtrees of $T$.
Thus, a $P$-element in a left subtree of $T_u$ is greater than a $P$-element in a right subtree of $T_u$ if comparable.
So, the above two requirements (i) and (ii) hold for this case.

We next examine the cases where $u$ is the root of a subtree of $T$.
Recall that the maximum $P$-element contained in a group marker
is the maximum of the whole group.
Then, the maximum $P$-element must be the root of the corresponding subtree of $T$
formed by the $P$-elements in the group in view of the building process of $T$.
Without loss of generality, we take the rightmost type I group, i.e., the one
with $C_k$ as the group marker, to continue the analysis.
In this case, $u=s_{km_k}$.
The requirement (i) is clear since $s_{k m_k}$ is the maximum $P$-element.
As for the requirement (ii),
suppose in the linear extension $C_1 C_2 \cdots C_k$, the chains contained in the
subsequence $C_l C_{l+1} \cdots C_k$ constitute the type I group with $C_k$
as the group marker.
Noticing that when restricted to this subsequence, the
constructed plane tree is exactly the subtree $T_u$ with $s_{k m_k}$ as the root.
Then the requirement (ii) concerning the vertices in the subtrees of $T_u$ follows
by the same token
as that for the subtrees of $T$ verified above.

Iterating the above reasoning for $u$ being a vertex in $T$ in a kind of ``top-down" manner,
we conclude that the requirements (i) and (ii) hold for all vertices.
Therefore, $e$, the reverse of the preorder of $T$ is a linear extension of $P$, completing
the proof.  \qed

Now it is not hard to piece all properties above together to arrive at Theorem~\ref{thm:nc-better}.
Obviously, when $P$ is itself a chain, all inequalities become equalities whence the sharpness claim.

\section{MHCD and automorphism group}

\begin{definition}
	Let $\mathcal{C}=\{C_1,C_2,\ldots, C_n\}$ be an HCD of $P$.
	The $\mathcal{C}$-graph of $P$ is the graph $G_{\mathcal{C}}$ having $C_i$'s as vertices where $C_i$ and $C_j$ are adjacent if they are comparable.
\end{definition}

We denote by
$G_P$ the $\mathcal{C}$-graph corresponding to the MHCD of $P$ hereafter.
The function $Min(P)$ (and many other functions) on a poset $P$ is Lipschitz continuous,
i.e., for any element $z\in P$, $|Min(P)-Min(P\setminus \{z\})|\leq 1$.
However, $Min_h(P)$ does not necessarily share this property.
In fact, we have the characterization below.

\begin{theorem}
	Let $P$ be a poset and $z\in P$.
	Then, we have the sharp bounds
	\begin{align}
		Min_h(P\setminus \{z\}) \leq Min_h(P) \leq 2 Min_h(P\setminus \{z\})+1.
	\end{align}
\end{theorem}
\proof
Let $\mathcal{C}=\{C_1,C_2,\ldots, C_k\}$ be the MHCD of $P$.
Deleting $z$ from the chain $C_i$ containing it, we obtain an HCD
$\mathcal{C'}=\{C_1,\ldots, C_i\setminus\{z\},\ldots, C_k\}$ of $P\setminus \{z\}$.
Thus, $Min_h(P\setminus \{z\}) \leq Min_h(P)$.
It is clear that if there is no edge in $G_{\mathcal{C}}$ and $|C_j|>1$ for all $1\leq j \leq k$,
$\mathcal{C'}$ will be the MHCD of $P\setminus \{z\}$ whence
the equality can be achieved.

For the second inequality, let $\mathcal{C'}$ and $\mathcal{C}$ be
the MHCDs of $P\setminus \{z\}$ and $P$, respectively.
We first claim that

\emph{Claim.} The elements of each chain in $\mathcal{C'}$ can be contained in at most
two chains in $\mathcal{C}$.\\
If otherwise, there exists a chain $C'_i$ in $\mathcal{C'}$ whose elements are contained in at least
three different chains $C_1,\: C_2,\: C_3$ in $\mathcal{C}$.
By the pigeonhole principle, there are at least two of them either comparable to $z$ or not.
W.l.o.g., suppose $C_1,\: C_2$ are both comparable to $z$. Then $C_1,\: C_2$ are both comparable to
the chain containing $z$.
Now, for any chain $C \in \mathcal{C}$ such that $z\notin C$, if $C_1$ is comparable to $C$, so is $C_2$, since
$C_1$ and $C_2$ are contained in the same chain $C'_i$, and vice versa.
Thus, $C_1$ and $C_2$ can be combined in $\mathcal{C}$ will induce a new HCD with less number of chains by merging $C_1$ and $C_2$,
which contradicts the assumption that
$\mathcal{C}$ is minimum. Hence, the claim is affirmed.

\begin{figure}[!ht]
	\centering
	\includegraphics[scale=.78]{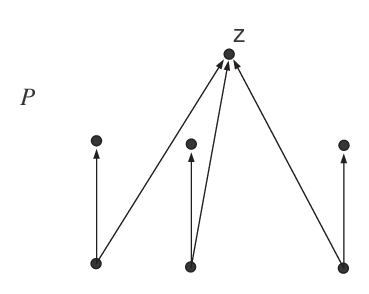}
	\caption{An example achieves the upper bound of $Min_h(P)$.}\label{4fig2}
\end{figure}

Note that $z$ may possibly form a chain itself in $\mathcal{C}$.
Therefore, we have
$$
Min_h(P) \leq 2 Min_h(P\setminus \{z\})+1.
$$
As far as the sharpness of the bound, consider the case:
suppose $\mathcal{C'}=\{C'_1,\ldots, C'_k\}$ with $|C'_i|=2$ for $1\leq i \leq k$, and there is no edge in
$G_{P\setminus \{z\}}$. Assume for each $1\leq i \leq k$, $z$ is larger than the minimum element
in $C_i$ but not comparable to the maximum element in $C_i$.
See Figure~\ref{4fig2} for an illustration.
Then, it is not difficult to check that each chain in the MHCD of $P$ contains only one element.
Thus, $Min_h(P)=2k+1=2 Min_h(P\setminus \{z\})+1$.
This completes the proof. \qed

Taking advantage of the uniqueness of the MHCD, we are able to
provide a characterization of poset automorphisms next.
See the papers~\cite{barmak,GB,BS,stan2} and references therein for previous studies
on poset automorphisms.
The following notation will be assumed through the end of the section.
Let $\Auto(P)$ denote the automorphism group of $P$, and let $\Auto(G_P)$ denote the automorphism group of the graph $G_P$. Let $\mathcal{C}=\{C_1,C_2,\dots, C_m\}$ be the MHCD of $P$.
Without loss of generality, assume the chains in the set $S_i=\{C_{k_i+1},C_{k_i+2},\dots, C_{k_{i+1}}\}$ have the same length $l_i$ for $0\leq i \leq t-1$, and $l_i\neq l_j$ if $i\neq j$, where $t\geq 1$, $k_0=0$, $k_t=m$, and $\sum_{i=0}^{t-1} l_i (k_{i+1}-k_{i}) =|P|$.
Denote by $O_P$ the group generated by all permutations $\pi$ on $V(G_P)$ such that every cycle of $\pi$ contains elements from the same set $S_i$ for some $0\leq i \leq t-1$.

\begin{theorem}\label{4thm2}
	There exists an acyclic orientation $\Asyc(G_P)$ of $G_P$ such that the automorphism group $\Auto(P)$ is isomorphic to a subgroup of the group
	$$
	\Auto(\Asyc(G_P))\bigcap O_P.
	$$
\end{theorem}

\proof
Note that both $\Auto(\Asyc(G_P))$ and $O_P$ are groups of permutations on $V(G_p)$.
Thus $\Auto(\Asyc(G_P))\bigcap O_P$ is indeed a group.
We first show that $\Auto(P)$ is isomorphic to a subgroup of $O_P$.

{\em Claim~$1$}. Under any $g\in \Auto(P)$, two $P$-elements in the same chain $C_i$ of the MHCD must be mapped to two $P$-elements in the same chain $C_j$ of the MHCD with $|C_i|=|C_j|$.\\
It is easy to see that any automorphism $g$ will map a chain to a chain, thus a chain decomposition to a chain decomposition with the same number of chains. The resulting chain decomposition $\mathcal{C'}$ induced by $g$ acting on the MHCD must be an HCD.
This is verified as follows.
Let $x_i'$ be an element from a chain $C_i'$ of $\mathcal{C'}$ and $x_j'$
be an element from a chain $C_j'$ of $\mathcal{C'}$. Suppose $x_i'$ and $x_j'$ are comparable.
Note that $g^{-1}\in \Auto(P)$.
Then their respective preimages $x_i=g^{-1}(x_i')$ and $x_j=g^{-1}(x_j')$ must be comparable as well.
Let $y_j'$ be any other element from $C_j'$ and let $y_j=g^{-1}(y_j')$.
Note that by construction $y_j$ is contained in the same chain as $x_j$.
So $x_i$ and $y_j$ are comparable. Thus their images under $g$, $x_i'$
and $y_j'$ must be comparable too, whence $C_i'$ and $C_j'$ are comparable.
Similarly, if $x_i'$ and $x_j'$ are incomparable, $C_i'$ and $C_j'$
are incomparable. Therefore, $\mathcal{C'}$ is an HCD.
By the uniqueness of the MHCD, each chain in the MHCD must be
uniquely and exclusively mapped to a chain in the MHCD under $g$, whence Claim~$1$.

Following Claim~$1$, any chain in $S_i$ will be mapped to a chain in $S_i$. Furthermore, the maximum element of a chain there must be mapped to the maximum element of a chain, and for each chain $C_i$ of the MHCD, once the image of the maximum element in the chain under $g$ is determined, then $g$ is completely determined. 
Thus, $g$ uniquely induces an element in $O_P$.
For $g,~g'\in \Auto(P)$, it is not hard to see that
the induced element of $g\circ g'$ is also in $O_P$.
Therefore, $\Auto(P)$ is isomorphic to a subgroup of $O_P$.

Let $\mathcal{C}=\{C_1,C_2,\ldots, C_n\}$ be the MHCD of a poset $P$.
Then, $(\mathcal{C},\leq_h)$ is obviously a well-defined poset, where the relation $\leq_h$ is defined as follows:
$C_i \leq_h C_i$, and for $i\neq j$, $C_i \leq_h C_j$ iff $\min(C_i)<\min(C_j)$.
Let $\Asyc(G_P)$ be the orientation induced by the poset $(\mathcal{C},\leq_h)$,
i.e., the edge between $C_i$ and $C_j$ in $G_P$ is oriented from $C_i$ to $C_j$
if $C_i\leq_h C_j$. We next show

{\em Claim~$2$}. $\Auto(P)$ is isomorphic to a subgroup of $\Auto(\Asyc(G_P))$.\\
To prove this, we show that each $g\in \Auto(P)$ uniquely induces an automorphism $\tilde{g}\in \Auto(\Asyc(G_P))$. From Claim~$1$, we know that $g$ maps chain to chain so that it induces a bijection $\tilde{g}$ on $V(G_P)$. It suffices to verify that $\tilde{g}$ preserves directed edges in $\Asyc(G_P)$.
Given a directed edge $C_i\rightarrow C_j$, by construction $\min(C_i)<\min(C_j)$. Thus,
$g(\min(C_i))<g(\min(C_j))$, and there is an edge between the chain $C'_i$ containing $g(\min(C_i))$ and the chain $C'_j$ containing $\min(C_j)$. By construction of the orientation,  $C'_i$ is directed to $C'_j$. Note under $\tilde{g}$, $C_i$ and $C_j$ will be mapped to $C'_i$ and $C'_j$ respectively. Obviously, $\tilde{g}$ maps non-adjacent pairs to non-adjacent pairs. 
Hence, $\tilde{g}\in \Auto(\Asyc(G_P))$, whence the claim.
Finally, for $g,~g'\in \Auto(P)$, it is easy to check that $g\circ g'$ induces 
$$
\tilde{g}\circ \tilde{g'}\in \Auto(\Asyc(G_P)) .
$$
Therefore, Claim~$2$ holds, and the theorem follows.\qed

Obviously, $\Auto(\Asyc(G_P))\subset \Auto(G_P)$, then the following corollary holds.
\begin{corollary}\label{mainthm}
	The automorphism group $\Auto(P)$ is isomorphic to a subgroup of the group
	$$
	\Auto(G_P)\bigcap O_P.
	$$
\end{corollary}

We end this paper with some {\bf future study problems}: (1) in which more general cases can some of the equalities be
achieved in Theorem~\ref{thm:nc-better}, e.g., when $Min_{nc}(P)=Min_h(P)$? (2) how many
noncrossing partitions are there for a given poset $P$? Note that the answer is given by the famous Catalan numbers
when $P=[n]$, and
(3) following the notation of Theorem~\ref{4thm2}, consider
the sufficient and/or necessary conditions such that
$$
\Auto(P) \sim \Auto(\Asyc(G_P))\bigcap O_P.
$$

\section*{Declarations}

{\bf Funding Statement:} This work was supported by the Anhui Provincial Natural Science Foundation of China (No.~2208085MA02)
and Overseas Returnee Support Project on Innovation and Entrepreneurship of Anhui Province (No.~11190-46252022001).

\noindent{\bf Conflict of Interest:} None.

\noindent{\bf Data Availability:} None.

\noindent{\bf Author Contribution:} R.X.F.C. did the research and wrote the paper.

	\section*{Acknowledgements}

The author would like to thank the anonymous referees for valuable comments and suggestions
which improved the presentation of the paper.

\end{document}